\documentclass[12pt,reqno,a4wide]{amsart}

\usepackage{tikz} 

\usetikzlibrary{decorations.pathreplacing}
\usetikzlibrary{arrows,decorations.pathmorphing,snakes}
\usetikzlibrary{fit}
\usepackage{calrsfs}
\usepackage[charter]{mathdesign}

\allowdisplaybreaks

\newtheorem{theorem}{Theorem}

\newtheorem{remark}{Remark}

\begin{document}
	
\title[On some problems    about ternary paths]{On some problems    about ternary paths --- a linear algebra approach  }

\author[Helmut Prodinger]{Helmut Prodinger}
\address{Department of Mathematics, University of Stellenbosch 7602,
	Stellenbosch, South Africa}
\email{hproding@sun.ac.za}


\begin{abstract} 
Ternary paths consist of an up-step of one unit, a down-step of two units, never go below the $x$-axis, and return to the $x$-axis.
This paper addresses the enumeration of partial ternary paths, ending at a given level $i$, reading the path either from left to right or from right to left. Since the paths are not symmetric w.r.t.\ left vs.\ right, as classical Dyck paths, this leads to different results. The right to left enumeration  is quite challenging, but leads at the end to very satisfying results. The methods are elementary (solving systems of linear equations). In this way, several conjectures left open in Naiomi Cameron's Ph.D. thesis could be successfully settled. 
\end{abstract}

\maketitle

\section{Introduction}

Ternary paths are cousins of Dyck paths, but with up-steps $(1,1)$ and down-steps $(1,-2)$, starting at the origin and never going below the $x$-axis. In most cases, one is interested in such paths that also end at the $x$-axis, but also at paths ending at level $i$
after $n$ steps.

Ternary paths are also called $2$-Dyck paths, see \cite{Selkirk-master}, where more general families of paths are studied.

Ternary paths are of interest at least for the following reasons:

\begin{itemize}
	\item They are no longer symmetric with respect to left $\leftrightarrow$ right, i.e., if the up-step is $(1,2)$ and the downstep is $(1,-1)$
	one gets slightly different results. In particular, arguments using symmetry as often employed when dealing with Dyck paths are no longer possible.
	
	\item Although some underlying generating functions are cubic, they are still manageable, due to a substitution, which allows to separate one factor, and dealing with a quadratic equation only. There is some recent interest in lattice paths models where the generating function is cubic, see e. g. \cite{ProSelWag}.
	
	\item Although one can look at more general classes of paths, the ternary case is of a nature to allow for explicit results.
	
	\item Knuth~\cite{christ} based his always popular christmas talk on the related concept of ternary trees, mentioning that he was and is interested in the subject for some 50 years.
	
	\item Ternary paths form a large portion of N.~Cameron's thesis~\cite{Cameron} but some answers were formulated in forms of conjectures.
\end{itemize}

Our method of choice is to find generating functions for ternary paths bounded by $h$ (thus the second coordinate is never $>h$) and then letting $h\to\infty$. This has the advantage that one has to deal only with a finite set of linear equations, and it can be solved explicitly using Cramer's rule. 

There is another (more high level) approach based on the \emph{kernel method}, and we might come back to that in a future publication. As one reviewer suggested, it would be interesting to see some of this here, but it would almost double the length of the paper, and would not stylistically fit. The interested reader will find a recent application of the kernel method to a ternary structure in \cite{kernel-Prodinger}, namely   the so-called S-Motzkin paths \cite{ProSelWag}, closely related to ternary paths, but only introduced recently.

We will address the following questions: Enumeration of ternary paths from left to right, starting at level 0 and ending at level $i$,
and enumeration of ternary paths from  right to left, starting at level 0 and ending at level $i$. The second question is harder than the first one. As we will see, from a certain cubic equation, the first root is responsible for the left to right enumeration, while the other two combined  are responsible for the right to left enumeration.

As a corollary, we will compute the (cumulated) area, summed over all ternary paths of length $3n$, again confirming a conjecture of N.~Cameron's. For a given path 
$(0,c_0), (1,c_1),\dots,\linebreak (3n,c_{3n})$, the area is defined to be the sum of the ordinates: $c_0+c_1+\cdots+c_{3n}$.

Banderier and Gittenberger study the area in a more general setting \cite{BG}, but the type of explicit results that we obtain   is restricted to the binary (Dyck)  base and the ternary case, as in this paper.

Here are our main findings. We only state results of an explicit enumeration type; perhaps it is more important to obtain the relevant generating functions, which will appear in later sections.

\begin{theorem} Enumeration of partial ternary paths and area.

	\begin{itemize}

		\item The number of ternary paths (from left to right) of length $3N+i$, ending on level $3K+i$, is given by
		\begin{equation*}
\binom{3N+i+1}{N-K}-3\binom{3N+i}{N-K-1}.
		\end{equation*} 
		\item The number of ternary paths (from  right to left) of length $3N+2i$, ending on level $i$, is given by
		\begin{equation*}
		\sum_{0\le k\le i/2}\binom{i-k}{k}\bigg[\binom{3N+2i+1}{N+k}-3\binom{3N+2i}{N+k-1}\bigg].
\end{equation*}		
\item The total area of all ternary paths of length $3N$, is given by
\begin{equation*}
\sum_{k=0}^{N-1}3^{k+1}\binom{3N-k}{N-1-k}.
\end{equation*}

	\end{itemize}
	
\end{theorem}

\textbf{Remarks.}
	\begin{itemize}
		
		\item [$\heartsuit$]
	Part 1 (left to right enumeration) was known to N.~Cameron~\cite{Cameron}, but, according to \cite[1.4.5]{Krattenthaler-survey}
	was redicovered by many people over the years.

		\item [$\diamondsuit$]
	Part 2 (right to left enumeration) was covered in \cite[1.4.7]{Krattenthaler-survey} in a different form. It would require a little bit of effort to show directly that the different forms are equivalent.

		\item [$\spadesuit$]
  The lattice path model in	 \cite[Section 4]{Krattenthaler-survey} is equivalent to ternary paths, but different: It consists of steps north and east of unit length each, together with some boundaries. The notion of \emph{area}, as of interest to N.~Cameron~\cite{Cameron} cannot be expressed in this model (at least not in an obvious way).

		\item[$\clubsuit$] 
Our approach is by first counting ternary paths with height restrictions. While this is only a vehicle here to get to the explicit generating functions, the results are of independent interest, and the explicit knowledge of the roots of the cubic equations involved, seems to be essential. For the more traditional class of Dyck paths, the average height (in terms of planted plane trees)	was treated in the seminal paper ~\cite{deBKR}.  Kemp~\cite{Kemp-prefix} considered the average height of prefixes of the Dyck-language (Dyck paths). To engage into a similar analysis of ternary paths or prefixes of ternary paths requires an explicit knowledge of the generating functions just mentioned. This might be a project of the future.
\end{itemize}

\section{A warm-up: the equation for ternary trees}

We start with the equation  $X=1+xX^3$ for ternary trees. 
The combinatorially relevant solution is $\mathcal{B}_3(x)$ and satisfies indeed 
 the equation  $\mathcal{B}_3(x)=1+x\mathcal{B}_3(x)^3$, since a ternary tree is either the empty tree or a root, followed
 by 3 (ordered) subtrees, which are also ternary trees. For $t$-ary trees, one has $\mathcal{B}_t(x)=1+x\mathcal{B}_t(x)^t$,
 but explicit results are only expected for $t=2$ and $t=3$. In fact, N.~Cameron's thesis \cite{Cameron} concentrates on binary and ternary trees and lattice paths as well.

In order to describe the other two solutions which are linked to what Knuth calls $(3/2)$-ary trees in \cite{christ}, we use the notation as given in \cite{GKP}:
\begin{equation*}
\mathcal{B}_t(x)^r=\sum_{k\ge0}\binom{tk+r}{k}\frac{r}{tk+r}x^k.
\end{equation*}
The combinatorially interesting solution is $\mathcal{B}_3(x)$. Using the substitution $x=t(1-t)^2$, the other two solutions are 
\begin{equation*}
\sigma_{1}=	-\frac1{2(1-t)}-\frac{\sqrt{4t-3t^2}}{2t(1-t)},\quad
\sigma_{2}=	-\frac1{2(1-t)}+\frac{\sqrt{4t-3t^2}}{2t(1-t)}.
\end{equation*}
The first part $\big(-\frac1{2(1-t)}=-\frac12\mathcal{B}_3(x)\big)$ is well understood, so we will consider $-\frac{\sqrt{4t-3t^2}}{2t(1-t)}$ 
and make use of
\begin{equation*}
t^k=\sum_{n\ge k}\binom{3n-k-1}{n-k}\frac knx^n,
\end{equation*}
which follows from the Lagrange inversion formula. Now we compute
\begin{align*}
-x^{1/2}\frac{\sqrt{4t-3t^2}}{2t(1-t)}&=-\frac{\sqrt{4t-3t^2}}{2t^{1/2}}=-\sqrt{1-\tfrac34t}=\sum_{k\ge0}(-1)^{k-1}(\tfrac34)^k\binom{1/2}{k}t^k\\
&=\sum_{k\ge0}(-1)^{k-1}(\tfrac34)^k\binom{1/2}{k}\sum_{n\ge k}\binom{3n-k-1}{n-k}\frac knx^n\\
&=\sum_{n\ge0}x^n	\sum_{1\le k\le n}(-1)^{k-1}(\tfrac34)^k\binom{1/2}{k}\binom{3n-k-1}{n-k}\frac kn\\
&=\sum_{n\ge0}x^n{\frac {{27}^{n}\Gamma  ( n+1/6 ) \Gamma  ( n-1/6 ) }{12\pi \Gamma  ( 2n+1 ) }}\\
&=-\sum_{n\ge0}x^n \binom{3n/2-1/2}{n}\frac{1}{3n-1}.
\end{align*}
The simplification of the inner sum can be done by a computer. Therefore
\begin{equation*}
\sigma_{1,2}=-\frac12\sum_{n\ge0}\binom{3n}{n}\frac{1}{2n+1}x^n\pm\sum_{n\ge0}x^{n-\frac12} \binom{3n/2-1/2}{n}\frac{1}{3n-1}.
\end{equation*}
A direct calculation shows that
\begin{equation*}
x^{1/2}\sigma_{1}= \sum_{n\ge0}\binom{3n/2-1/2}{n}\frac{1}{3n-1}x^{n/2}(-1)^n=-\mathcal{B}_{3/2}(-x^{1/2})^{-1/2}
\end{equation*}
and
\begin{equation*}
x^{1/2}\sigma_{2}= -\sum_{n\ge0}\binom{3n/2-1/2}{n}\frac{1}{3n-1}x^{n/2}=\mathcal{B}_{3/2}(x^{1/2})^{-1/2}.
\end{equation*}
Finally, the equation of interest factors as
\begin{equation*}
x(X-\mathcal{B}_3(x))(X-\sigma_1)(X-\sigma_2)=xX^3+1-X.
\end{equation*}
Thus $x\mathcal{B}_3(x)\sigma_1\sigma_2=-1$,
which leads to
\begin{align*}
\mathcal{B}_3(x)=\Big(\mathcal{B}_{3/2}(x^{1/2})^{-1/2}\mathcal{B}_{3/2}(-x^{1/2})^{-1/2}\Bigr)^{-1}=
\mathcal{B}_{3/2}(x^{1/2})^{1/2}\mathcal{B}_{3/2}(-x^{1/2})^{1/2}.
\end{align*}

\begin{remark}
	This factorization had been obtained by Bousquet-M\'elou and Petkovsek~\cite{BM-P} using different methods. The  approach   of choice in the present paper is to use the substitution $x=t(1-t)^2$ which is central for the later sections. 
\end{remark}

\begin{remark}
S. Selkirk in \cite{Selkirk-master} was able to factor the more general equation $X=1+xX^d$.
	\end{remark}

\section{Enumeration of ternary paths from left to right}

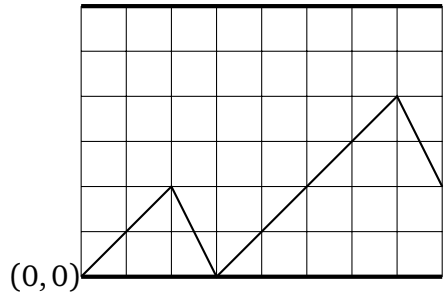
\begin{center}
	\begin{tikzpicture}[scale=0.6]
	
	\draw[step=1.cm,black] (-0.0,-0.0) grid (8.0,6.0);

	\draw[ultra thick] (0,0) to (8,0);
	\draw[ultra thick] (0,6) to (8,6);
	
	\draw[thick] (0,0) -- (1,1) -- (2,2)-- (3,0)-- (4,1)-- (5,2)-- (6,3)-- (7,4)-- (8,2) ;
	
	\node at (-0.8,0.0){$(0,0)$};
	
	\end{tikzpicture}
	\end	{center}

The picture shows a ternary path ending in $(8,2)$ and being bounded by $6$.

Let $a_{n,k}$ be the number of ternary paths ending at $(n,k)$ and being bounded by $h$. 
In order not to clutter the notation, we did not put the letter $h$ into the definition, especially, 
since $h$ has only a very temporary meaning here. 

The recursion (for $n\ge1$)
\begin{equation*}
a_{n,k}=a_{n-1,k-1}+a_{n-1,k+2}
\end{equation*}
with the understanding that $a_{n,k}$ should be interpreted as 0 if $k<0$ or $k>h$ is clear, distinguishing the two instances that can lead to level $k$ in one step. The starting value
is $a_{0,0}=1$. It is natural to introduce the generating functions
\begin{equation*}
f_k=f_k(z)=\sum_{n\ge0}a_{n,k}z^n,
\end{equation*}
then
$f_k=zf_{k-1}+zf_{k+2}+[\!\![k=0]\!\!]$,\footnote{The Iverson notation $[\!\![P]\!\!]$, very common these days, is 1 if condition $P$ is true and 0 if condition $P$ is false, see \cite{GKP}.} which is best written as a matrix equation

\begin{equation*}
\begin{bmatrix}
1&0& -z&\dots\\
-z&1&0& -z&\dots\\
0&-z&1&0& -z&\dots\\
&&&\vdots\\
&&&& -z&1
\end{bmatrix}
\begin{bmatrix}
f_0\\f_1\\f_2\\ \vdots\\f_h
\end{bmatrix}
=\begin{bmatrix}
1\\0\\0\\ \vdots\\0
\end{bmatrix}
\end{equation*}
Here is a little list (with $h$ being chosen large enough).
\begin{align*}
	f_0(z)&=1+{z}^{3}+3 {z}^{6}+12 {z}^{9}+55 {z}^{12}+273 {z}^{15}+\cdots,\\
	f_1(z)&=z+2 {z}^{4}+7 {z}^{7}+30 {z}^{10}+143 {z}^{13}+\cdots,\\
	f_2(z)&={z}^{2}+3 {z}^{5}+12 {z}^{8}+55 {z}^{11}+273 {z}^{14}+\cdots,\\
	f_3(z)&={z}^{3}+4 {z}^{6}+18 {z}^{9}+88 {z}^{12}+455 {z}^{15}+\cdots,\\
	f_4(z)&={z}^{4}+5 {z}^{7}+25 {z}^{10}+130 {z}^{13}+\cdots,\\
	f_5(z)&={z}^{5}+6 {z}^{8}+33 {z}^{11}+182 {z}^{14}+\cdots,\\
	f_6(z)&={z}^{6}+7 {z}^{9}+42 {z}^{12}+245 {z}^{15}+\cdots\\
\end{align*}

Now let $d_h$ be the determinant of this matrix with $h$ rows and columns. We have $d_0=1$, $d_1=1$, $d_2=1$, and the recursion
\begin{equation*}
d_h=d_{h-1}-z^3d_{h-3},
\end{equation*}
which is obtained by expanding the determinant along the first column, say. The method is described in \cite{Prodinger-handbook}.
The characteristic equation of this recursion is the cubic equation
\begin{equation*}
\lambda^3-\lambda^2+z^3=0.
\end{equation*}
Note also the generating function
\begin{align*}
R(X)&=\sum_{h\ge0}d_hX^h=1+X+X^2+\sum_{h\ge3}(d_{h-1}-z^3d_{h-3})X^h\\
&=1+X+X^2+X\sum_{h\ge2}d_{h}X^h-z^3X^3R(X)\\
&=1+X+X^2+XR(X) -X-X^2-z^3X^3R(X),
\end{align*}
or
\begin{equation*}
R(X)=\frac{1}{1-X+z^3X^3}=\sum_{j\ge0}d_{j}X^j.
\end{equation*}

This cubic equation becomes manageable with the substition $z^3=t(1-t)^2$, which featured prominently in \cite{ProSelWag}.
It has the 3 roots
\begin{align*}
r_1&=1-t,\quad r_{2}=\frac{t+\sqrt{4t-3t^2}}{2},\quad r_{3}=\frac{t-\sqrt{4t-3t^2}}{2},
\end{align*}
and now the relation to the cubic equation that we studied as a warm-up becomes apparent.

Cramer's rule now leads to
\begin{equation*}
f_k= z^{k}\frac{d_{h-k}}{d_{h+1}},
\end{equation*}
which, when performing the limit $h\to\infty$, leads to
\begin{equation*}
f_k(z)= z^{k}r_1^{-k-1}=\frac{z^{k}}{(1-t)^{k+1}}.
\end{equation*}
This form will be useful later, but we would also like to compute $[z^n]f_k(z)$, i. e., then numbers $a_{n,k}$. We can only have contributions (which is also clear for combinatorial reasons) if $n\equiv k \bmod 3$. So let us set $n=3N+i$, $k=3K+i$ for 
$i=0,1,2$, and compute
\begin{align*}
[z^{3N+i}]\frac{ z^{3K+i}}{(1-t)^{3K+i+1}}&=[z^{3N-3K}]\frac{1}{(1-t)^{3K+i+1}}.
\end{align*}
In order not to confuse matters, it helps to set $z^3=x=t(1-t)^2$. Then we can continue
\begin{align*}
a_{3N+i,3K+i}&=[x^{N-K}]\frac{1}{(1-t)^{3K+i+1}}\\
&=\frac1{2\pi i}\oint\frac{dx}{x^{N-K+1}}\frac{1}{(1-t)^{3K+i+1}}\\
&=\frac1{2\pi i}\oint\frac{dt(1-t)(1-3t)}{t^{N-K+1}(1-t)^{2N-2K+2}}\frac{1}{(1-t)^{3K+i+1}}\\
&=[t^{N-K}](1-3t)\frac{1}{(1-t)^{2N+K+i+2}}\\
&=\binom{3N+i+1}{N-K}-3\binom{3N+i}{N-K-1}.
\end{align*}
These coefficients can also be obtained using the Lagrange inversion formula; we find the approach using a contour integral, as seen in \cite{deBKR}, more versatile. Some of the shortest proofs of the Lagrange inversion formula use indeed contour integration, see
\cite{Prodinger-handbook}.

Notice in particular the enumeration of paths ending at the $x$-axis:
\begin{align*}
a_{3N,0}
&=\binom{3N+1}{N}-3\binom{3N}{N-1}=\frac1{2N+1}\binom{3N}{N},
\end{align*}
a generalized Catalan number.

We can check now that
\begin{equation*}
f_{3K+i}(z)=\sum_{N\ge0}a_{3N+i,3K+i}z^{3N+i}
\end{equation*}
coincides with the previous list. 

\section{Enumeration of ternary paths from right to left}

While the left-to-right enumeration was done successfully in \cite{Cameron} (with a different approach), the enumeration from right to left was not as complete as what we are going to do now and partially in conjectural state.

We still prefer to work from left to right, so we change our setting as follows:

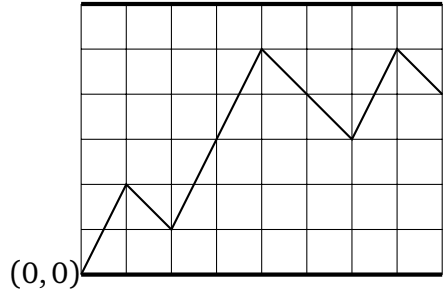
\begin{center}
	\begin{tikzpicture}[scale=0.6]
	
	\draw[step=1.cm,black] (-0.0,-0.0) grid (8.0,6.0);

	\draw[ultra thick] (0,0) to (8,0);
	\draw[ultra thick] (0,6) to (8,6);
	
	\draw[thick] (0,0) -- (1,2) -- (2,1)-- (3,3)-- (4,5)-- (5,4)-- (6,3)-- (7,5)-- (8,4) ;
	
	\node at (-0.8,0.0){$(0,0)$};
	
	\end{tikzpicture}
	\end	{center}

	The picture shows a reversed ternary path ending in $(8,4)$ and being bounded by $6$.

For the notation, we switch from $a_{n,k}$ to $b_{n,k}$ and from $f_k(z)$ to $g_k(z)$.

Here is a little list (the boundary $h$ is assumed to be large enough):
\begin{align*}
	g_0(z)&=1+{z}^{3}+3 {z}^{6}+12 {z}^{9}+55 {z}^{12}+273 {z}^{15}+\cdots,\\
	g_1(z)&={z}^{2}+3 {z}^{5}+12 {z}^{8}+55 {z}^{11}+273 {z}^{14}+1428 {z}^{
		17}+\cdots,\\
	g_2(z)&=z+3 {z}^{4}+12 {z}^{7}+55 {z}^{10}+273 {z}^{13}+1428 {z}^{16}+\cdots,\\
	g_3(z)&=2 {z}^{3}+9 {z}^{6}+43 {z}^{9}+218 {z}^{12}+1155 {z}^{15}+\cdots,\\
	g_4(z)&={z}^{2}+6 {z}^{5}+31 {z}^{8}+163 {z}^{11}+882 {z}^{14}+4896 {z}^{
		17}+\cdots,\\
	g_5(z)&=3 {z}^{4}+19 {z}^{7}+108 {z}^{10}+609 {z}^{13}+3468 {z}^{16}+\cdots.\\
	g_6(z)&={z}^{3}+10 {z}^{6}+65 {z}^{9}+391 {z}^{12}+2313 {z}^{15}+\cdots.
\end{align*}

The linear system changes now as follows:

\begin{equation*}
\begin{bmatrix}
1& -z&\dots\\
0&1& -z&\dots\\
-z&0&1& -z&\dots\\
&&&\vdots\\
&&&& 0&1
\end{bmatrix}
\begin{bmatrix}
g_0\\g_1\\g_2\\ \vdots\\g_h
\end{bmatrix}
=\begin{bmatrix}
1\\0\\0\\ \vdots\\0
\end{bmatrix}
\end{equation*}

The determinant of the matrix is the same as before by transposition: $d_{h+1}$. However, the application of Cramer's rule is more involved now. We must evaluate the determinant of
\begin{align*}
&\begin{bmatrix}
0& 0&\dots& 1 &\dots&0\\
0&1& -z&0&\dots\\
-z&0&1& 0&\dots\\
&&&\vdots\\
&&&0& 0&1
\end{bmatrix}\\
&\underbrace{\phantom{xxxxxxxxx}}\quad\underbrace{\phantom{xxxxxx}}\\
&\hspace*{1.2cm}i\hspace*{2.2cm}j
\end{align*}
with $i+j=h$. Call this determinant $(-1)^i\Delta_{i,j}$. We want to find a recursion for it.

By expansion, we find the recursion 
\begin{equation*}
\Delta_{i,j}=z\Delta_{i-2,j}-z^3\Delta_{i-3,j}.
\end{equation*}
The characteristic equation of this recursion is
\begin{equation*}
Y^3-zY+z^3=0.
\end{equation*}
Setting $Y=\frac{z^2}{X}$, this leads to
\begin{equation*}
X^3-X^2+z^3=0,
\end{equation*} 
which was the equation studied before. 

The expansion from the other end leads to 
\begin{equation*}
\Delta_{i,j}=\Delta_{i,j-1}-z^3\Delta_{i,j-3}.
\end{equation*}

In particular, 
\begin{equation*}
	\Delta_{0,j}=d_j,\quad \Delta_{1,j}=-z^2d_{j-1}\ (j\ge1),\quad \Delta_{1,0}=0,\quad \Delta_{2,j}=zd_j.
\end{equation*}

Now we compute 
\begin{align*}
F(X,Y)&:=\sum_{i\ge0}\sum_{j\ge0}\Delta_{i,j}X^iY^j\\*
&=\sum_{j\ge0}\Delta_{0,j}Y^j
+X\sum_{j\ge0}\Delta_{1,j}Y^j
+X^2\sum_{j\ge0}\Delta_{2,j}Y^j
+\sum_{i\ge3}\sum_{j\ge0}\Delta_{i,j}X^iY^j\\
&=\sum_{j\ge0}\Delta_{0,j}Y^j
-z^2X\sum_{j\ge1}\Delta_{0,j-1}Y^j
+zX^2\sum_{j\ge0}\Delta_{0,j}Y^j
+\sum_{i\ge3}\sum_{j\ge0}[z\Delta_{i-2,j}-z^3\Delta_{i-3,j}]X^iY^j\\
&=R(Y)-z^2XYR(Y)+zX^2R(Y)+z\sum_{i\ge3}\sum_{j\ge0}\Delta_{i-2,j}X^iY^j\\
&-z^3\sum_{i\ge3}\sum_{j\ge0}\Delta_{i-3,j}X^iY^j\\
&=R(Y)(1-z^2XY+zX^2)+zX^2\sum_{i\ge1}\sum_{j\ge0}\Delta_{i,j}X^iY^j-z^3X^3\sum_{i\ge0}\sum_{j\ge0}\Delta_{i,j}X^iY^j\\
&=R(Y)(1-z^2XY+zX^2)-zX^2\sum_{j\ge0}\Delta_{0,j}Y^j-z^3X^3F(X,Y)\\
\end{align*}
whence
\begin{equation*}
F(X,Y)=\frac{1-z^2XY}{1-Y+z^3Y^3}\frac{1}{1-zX^2+z^3X^3}.
\end{equation*}
We can now continue with the computation, according to Cramer's rule:
\begin{align*}
(-1)^ig_i&=\frac1{d_{h+1}}[X^iY^{h-i}]F(X,Y)\\
&=\frac1{d_{h+1}}[X^iY^{h-i}]\frac{1-z^2XY}{1-Y+z^3Y^3}\frac{1}{1-zX^2+z^3X^3}\\
&=\frac1{d_{h+1}}[X^i]\frac{1}{1-zX^2+z^3X^3}[Y^{h-i}]\frac{1-z^2XY}{1-Y+z^3Y^3}\\
&=\frac1{d_{h+1}}[X^i]\frac{1}{1-zX^2+z^3X^3}[d_{h-i}-z^2Xd_{h-i-1}]\\
&=\frac{d_{h-i}}{d_{h+1}}[X^i]\frac{1}{1-zX^2+z^3X^3}-z^2\frac{d_{h-i-1}}{d_{h+1}}[X^{i-1}]\frac{1}{1-zX^2+z^3X^3}.
\end{align*}

Writing
\begin{equation*}
\frac{1}{1-zX^2+z^3X^3}=\sum_{n\ge0}w_nX^n,
\end{equation*}
we see that the numbers satisfy the recursion
\begin{equation*}
w_n-zw_{n-2}+z^3w_{n-3}=0.
\end{equation*}
The characteristic equation of this recursion is
\begin{equation*}
\lambda^3-z\lambda+z^3=0.
\end{equation*}
Setting $\lambda=z^2\mu$, this becomes
\begin{equation*}
1-\mu+z^3\mu^3=0
\end{equation*}
with some nice roots 
\begin{align*}
\mu_1&=\frac1{1-t},\quad \mu_2=\frac{t+\sqrt{4t-3t^2}}{2t(t-1)},\quad \mu_3=\frac{t-\sqrt{4t-3t^2}}{2t(t-1)}.
\end{align*}

Thus, in terms of $\lambda$, the roots are 
$z^2\mu_1$, $z^2\mu_2$, $z^2\mu_3$.
Consequently
\begin{equation*}
	w_n=\alpha(z^2\mu_1)^n+\beta(z^2\mu_2)^n+\gamma(z^2\mu_3)^n=z^{2n}(\alpha \mu_1^n+\beta \mu_2^n+\gamma \mu_3^n)
\end{equation*}
with
\begin{align*}
\alpha&=\frac{t}{3t-1},\\
\beta&=\frac{(2t-1)(3t-4)-\sqrt{4t-3t^2}}{2(3t-1)(3t-4)},\\
\gamma&=\frac{(2t-1)(3t-4)+\sqrt{4t-3t^2}}{2(3t-1)(3t-4)}.
\end{align*}

Note also that
\begin{equation*}
\mu_2+\mu_3=\frac{1}{t-1}\quad\text{and}\quad\mu_2\mu_3=\frac{1}{t(t-1)}.
\end{equation*}

\section{Pushing the boundary $h$ to infinity}

We start from the formula
\begin{align*}
(-1)^i	g_i	&=\frac{d_{h-i}}{d_{h+1}}w_i-z^2\frac{d_{h-i-1}}{d_{h+1}}w_{i-1}\\
\end{align*}
enumerating ternary path (from right to left), but written in a reversed way, ending at level~$i$.

Now we push the boundary $h$ to infinity, i. e. we have no more horizontal boundary. We write again $g_i=g_i(z)$.

From the explicit formula for the determinants $d_h$, we can conclude that only one of the 3 roots survives, and
\begin{equation*}
\lim_{h\to\infty}\frac{d_{h-i}}{d_{h+1}}=\frac1{(1-t)^{i+1}}\quad\text{and}\quad
\lim_{h\to\infty}\frac{d_{h-i-1}}{d_{h+1}}=\frac1{(1-t)^{i+2}}.
\end{equation*}
So, after taking limits,
\begin{equation*}
(-1)^ig_i=\frac{1}{(1-t)^{i+1}}z^{2i}(\alpha \mu_1^i+\beta \mu_2^i+\gamma \mu_3^i)-
\frac{1}{(1-t)^{i+2}}z^{2i}(\alpha \mu_1^{i-1}+\beta \mu_2^{i-1}+\gamma \mu_3^{i-1}).
\end{equation*}
Since $\mu_1=\frac1{1-t}$, there is cancellation:
\begin{equation*}
	(-1)^ig_i=\frac{1}{(1-t)^{i+1}}z^{2i}(\beta \mu_2^i+\gamma \mu_3^i)-
	\frac{1}{(1-t)^{i+2}}z^{2i}(\beta \mu_2^{i-1}+\gamma \mu_3^{i-1}).
\end{equation*}
This can be simplified:
\begin{align*}
g_i=\frac{(-1)^iz^{2i}}{(1-t)^{i+1}(3t-4)}\Bigl[(t-2)(\mu_2^i+\mu_3^i)+(\mu_2^{i-1}+\mu_3^{i-1})\Bigr].
\end{align*}
This formula works for $i=0$ as well.

In the next section we will explain why the square bracket always contains the factor $3t-4$.

\section{Explicit expressions for the generating functions $g_i$}

Let
\begin{equation*}
 p_m=	\mu_2^m+\mu_3^m
\end{equation*}
and
\begin{equation*}
u_m=(t-2)p_{m}+p_{m-1}.
\end{equation*}
We have
\begin{equation*}
\frac{u_m}{3t-4}=\frac{\mu_2^{m+1}-\mu_3^{m+1}}{\mu_2-\mu_3},
\end{equation*}
which, once it is known, is easy to prove by induction.

We need one form of the Girard-Waring formula, see e. g. \cite{Gould}:
\begin{equation*}
\frac{X^m-Y^m}{X-Y}=\sum_{0\le k\le (m-1)/2}(-1)^k\binom{m-1-k}{k}(XY)^{k}(X+Y)^{m-1-2k}.
\end{equation*}
In our application:
\begin{equation*}
\frac{\mu_2^{m+1}-\mu_3^{m+1}}{\mu_2-\mu_3}=\sum_{0\le k\le m/2}(-1)^k\binom{m-k}{k}\frac1{t^k(t-1)^{m-k}}.
\end{equation*}
Summarizing
\begin{align*}
\frac{u_m}{3t-4}=\sum_{0\le k\le m/2}(-1)^k\binom{m-k}{k}\frac1{t^k(t-1)^{m-k}}
\end{align*}
and
\begin{equation*}
g_i=\frac{(-1)^iz^{2i}}{(1-t)^{i+1}}\sum_{0\le k\le i/2}(-1)^k\binom{i-k}{k}\frac1{t^k(t-1)^{i-k}}.
\end{equation*}

It is easy to check that $g_i(z)$ has only nonzero coefficients for $z^n$ when $n+i\equiv 0\mod3$.

The enumeration of $[z^{3N+2i}]g_i(z)$ can now be done:
\begin{align*}
[z^{3N+2i}]g_i(z)&=[z^{3N}]\frac{(-1)^i}{(1-t)^{i+1}}\sum_{0\le k\le i/2}(-1)^k\binom{i-k}{k}\frac1{t^k(t-1)^{i-k}}\\
&=[x^{N}]\sum_{0\le k\le i/2}(-1)^{k+1}\binom{i-k}{k}\frac1{t^k(t-1)^{2i+1-k}}\\
&=\sum_{0\le k\le i/2}(-1)^{k+1}\binom{i-k}{k}\cdot \frac1{2\pi i}\oint\frac{dx}{x^{N+1}}\frac1{t^k(t-1)^{2i+1-k}}\\
&=\sum_{0\le k\le i/2}(-1)^{k+1}\binom{i-k}{k}\cdot \frac1{2\pi i}\oint\frac{dt(t-1)(3t-1)}{t^{N+1}(1-t)^{2N+2}}\frac1{t^k(t-1)^{2i+1-k}}\\
&=\sum_{0\le k\le i/2}\binom{i-k}{k}[t^{N+k}]\frac{(1-3t)}{(1-t)^{2N+2i+2-k}}\\
&=\sum_{0\le k\le i/2}\binom{i-k}{k}\bigg[\binom{3N+2i+1}{N+k}-3\binom{3N+2i}{N+k-1}\bigg].
\end{align*}

\section{The area }

As a corollary of our previous findings, we consider the total area of ternary paths, summed over all such paths of the same length.
Recall that the area of a path is the sum of its ordinates. This settles a conjectural result in N.~Cameron's thesis \cite{Cameron} as well.

Each contribution $c_i$ to the area of a path $(0,c_0=0),\dots,(3n,c_{3n}=0)$ can be seen as splitting the ternary path into a path of
length $k$ (left to right) ending at level $i$ and a path of length $3n-k$ (right to left) also ending at level $i$.
Since we are working with generating functions, all possible such splittings are taken into account when taking the product of two such generating functions.

The cumulated area is thus given as (write again $z^3=x=t(1-t)^2$)
\begin{align*}
\textsc{area}&=[z^{3n}]\sum_{i\ge0}if_i(z)g_i(z)\\
&=[z^{3n}]\sum_{i\ge0}i  \frac{z^i}{(1-t)^{i+1}}\frac{(-1)^iz^{2i}}{(1-t)^{i+1}(3t-4)}\Bigl[(t-2)(\mu_2^i+\mu_3^i)+(\mu_2^{i-1}+\mu_3^{i-1})\Bigr]\\
&=[z^{3n}]\frac{3t}{(1-3t)^2(1-t)}\\
&=\frac1{2\pi i}\oint\frac{dx}{x^{n+1}}\frac{3t}{(1-3t)^2(1-t)}\\
&=\frac1{2\pi i}\oint\frac{dt(1-3t)(1-t)}{t^{n+1}(1-t)^{2n+2}}\frac{3t}{(1-3t)^2(1-t)}\\
&=\frac1{2\pi i}\oint\frac{dt}{t^{n}(1-t)^{2n+2}}\frac{3}{(1-3t)}\\
&=[t^{n-1}]\frac{1}{(1-t)^{2n+2}}\sum_{k\ge0}3^{k+1}t^k\\
&=\sum_{k=0}^{n-1}3^{k+1}\binom{3n-k}{n-1-k}.
\end{align*}
This is a proof of Conjecture 3.1 in N.~Cameron's thesis \cite{Cameron}.

\textbf{Acknowledgment.} I would like to thank S. Selkirk and S. Wagner for related although different discussions.

Thanks are due to C. Krattenthaler for valuable information about lattice path enumeration.
\bibliographystyle{plain}


\begin{thebibliography}{10}
	
	\bibitem{BG}
	C.~Banderier and B.~Gittenberger.
	\newblock Analytic combinatorics of lattice paths: enumeration and asymptotics
	for the area.
	\newblock In {\em Fourth {C}olloquium on {M}athematics and {C}omputer {S}cience
		{A}lgorithms, {T}rees, {C}ombinatorics and {P}robabilities}, Discrete Math.
	Theor. Comput. Sci. Proc., AG, pages 345--355. 2006.
	
	\bibitem{BM-P}
	M.~Bousquet-M\'{e}lou and M.~Petkov\v{s}ek.
	\newblock Walks confined in a quadrant are not always {D}-finite.
	\newblock {\em Theoret. Comput. Sci.}, 307(2):257--276, 2003.
	\newblock in: Random generation of combinatorial objects and bijective
	combinatorics.
	
	\bibitem{Cameron}
	N.~Cameron.
	\newblock {\em PhD Thesis: Random Walks, Trees and Extensions of {R}iordan
		Group Techniques}.
	\newblock Howard University, 2002.
	
	\bibitem{deBKR}
	N.~G. de~Bruijn, D.~E. Knuth, and S.~O. Rice.
	\newblock The average height of planted plane trees.
	\newblock In {\em Graph theory and computing}, pages 15--22. 1972.
	
	\bibitem{Gould}
	H.~Gould.
	\newblock The {G}irard-{W}aring power sum formulas for symmetric functions and
	{F}ibonacci sequences.
	\newblock {\em The Fibonacci Quarterly}, 37:135--140, 1999.
	
	\bibitem{GKP}
	R.~L. Graham, D.~E. Knuth, and O.~Patashnik.
	\newblock {\em Concrete Mathematics}.
	\newblock Addison-Wesley, Reading, MA, 1999.
	
	\bibitem{Kemp-prefix}
	R.~Kemp.
	\newblock On the average depth of a prefix of a {D}ycklanguage {$D_{1}$}.
	\newblock {\em Discrete Math.}, 36(2):155--170, 1981.
	
	\bibitem{christ}
	D.~Knuth.
	\newblock {D}onald {K}nuth's 20th {A}nnual {C}hristmas {T}ree {L}ecture:
	(3/2)-ary {T}rees.
	\newblock \textsf{https://youtu.be/P4AaGQIo0HY}, Dec 2014.
	
	\bibitem{Krattenthaler-survey}
	Christian Krattenthaler.
	\newblock Lattice path enumeration.
	\newblock In {\em Handbook of enumerative combinatorics}, Discrete Math. Appl.
	(Boca Raton), pages 589--678. CRC Press, Boca Raton, FL, 2015.

	\bibitem{Prodinger-handbook}
	Helmut Prodinger.
	\newblock Analytic methods.
	\newblock In {\em Handbook of enumerative combinatorics}, Discrete Math. Appl.
	(Boca Raton), pages 173--252. CRC Press, Boca Raton, FL, 2015.


\bibitem{kernel-Prodinger}
Helmut Prodinger.
\newblock
{ Enumeration of S-Motzkin paths from left to right and from right to left --- a kernel method approach}, preprint.
	
	\bibitem{ProSelWag}
	H.~Prodinger, S.~J. Selkirk, and S.~Wagner.
	\newblock On two subclasses of {M}otzkin paths and their relation to ternary
	trees.
	\newblock In V.~Pillwein \and C.~Schneider, editor, {\em Algorithmic
		Combinatorics -- Enumerative Combinatorics, Special Functions and Computer
		Algebra}. Springer, Austria, 2019.
	
	\bibitem{Selkirk-master}
	S.~Selkirk.
	\newblock {\em MSc-thesis: {O}n a generalisation of $k$-{D}yck paths}.
	\newblock Stellenbosch University, 2019.
	
\end{thebibliography}

\end{document}